\documentclass[letterpaper, 10 pt, conference]{ieeeconf}  

\IEEEoverridecommandlockouts                              

\let\labelindent\relax

\usepackage[dvipsnames]{xcolor}

\newcommand{\vo}[1]{\boldsymbol{#1}}
\newcommand{\omegab}{{\vo{\omega}}}

\usepackage{graphicx,cite} 
\usepackage{amsmath, enumitem,amssymb} 
\usepackage{empheq} 
\title{\LARGE \bf
Nonlinear Attitude Estimation for Small UAVs with Low Power Microprocessors}

\author{Sunsoo Kim$^{1}$ and Vaishnav Tadiparthi$^{2}$ and Raktim Bhattacharya$^{3}$
\thanks{$^{1}$Sunsoo Kim is a Ph.D student in the Department of Electrical and Computer Engineering,  Texas A\&M University, College Station, TX 77840, USA. Email: {\tt\small kimsunsoo@tamu.edu}}%
\thanks{$^{2}$Vaishnav Tadiparthi is a Ph.D student in the Department of Aerospace Engineering, Texas A\&M University, College Station, TX 77840, USA. Email:
        {\tt\small vaishnavtv@tamu.edu}}%
\thanks{$^{3}$Raktim Bhattacharya is with the Faculty of Aerospace Engineering, Texas A\&M University, College Station, TX 77840, USA. Email:
        {\tt\small raktim@tamu.edu}}%
}

\begin{document}
\maketitle
\thispagestyle{empty}
\pagestyle{empty}

\begin{abstract}
Among algorithms used for sensor fusion for attitude estimation in unmanned aerial vehicles, the Extended Kalman Filter (EKF) is the most commonly used for estimation. In this paper, we propose a new version of  $\mathcal{H}_2$ estimation called extended $\mathcal{H}_2$ estimation that can overcome the limitations of the extended Kalman Filter, specifically with respect to computational speed, memory usage, and root mean squared error. We formulate a new attitude-estimation algorithm, where the filter gain is designed offline about a nominal operating point, but the filter dynamics is implemented using the nonlinear system dynamics. We refer to this implementation of the $\mathcal{H}_2$ optimal estimator as the extended  $\mathcal{H}_2$ estimator.  The solution presented is tested on two cases, corresponding to slow and rapid motions, and compared against the EKF in the performance metrics mentioned above.
\end{abstract}

\section{INTRODUCTION}
As unmanned air vehicles keep getting smaller and cheaper, the need for computationally efficient attitude estimation is growing rapidly. Attitude estimation is critical for a component of the flight control system controlling these systems \cite{eure2013application, gebre2004design, kada2016uav, weibel2015small}. While states pertaining to translational motion can be easily recovered from sensor data, orientation of these vehicles cannot be directly obtained from the same.  In that context, sensor fusing algorithms are typically employed to estimate the attitude/orientation of these vehicles.

In the recent past, MEMS (Micro-Electro Mechanical Systems) sensors like MARG (Magnetic, Angular Rate, and Gravity) sensor and IMU (Inertial Measurement Unit) have become increasingly common because of their low costs and small sizes.
These combine to from a three-axis gyroscope, a three-axis accelerometer, and a three axis magnetometer.
The gyroscope measures the angular velocity of the vehicle, the accelerometer measures the acceleration of the vehicle, and the magnetometer measures the magnetic vector.  It is important to note that measurements from MEMS sensor are corrupted by noise and bias. Additionally, rapid movements and magnetic disturbances can temporarily influence the attitude calculations \cite{chang2008integrated,fan2018magnetic}.

Each sensor in MARG and IMU can independently estimate all elements of a vehicle's attitude without external signals.
This is observed most popularly in conventional navigation systems, where the attitude is calculated by integrating the angular velocity obtained from the gyroscope outputs. However, this is not accurate because the gyroscope has a bias, which results in the accumulation of attitude error over long periods of time \cite{geiger2008mems}. Tilt angle can be estimated by an accelerometer computing gravitational force and a magnetometer measuring magnetic field.  Note however that this result is affected by rapid accelerations and magnetic disturbance \cite{luinge2005measuring, frosio2008autocalibration} as well. Thus, to obtain a reliably accurate estimate of the vehicle's attitude, measurements from all the sensors are fused in a filtering framework.

Since attitude estimation is inherently a nonlinear estimation problem, a number of nonlinear sensor fusing algorithms have been proposed for attitude estimation \cite{crassidis2007survey}. Among these, filters which require high computation, like unscented filter and particle filter etc., will be ignored for the purpose of this paper, since the focus is on estimating attitude in low-cost embedded processors. Among the rest, one of the most widely used techniques is the Extended Kalman filter (EKF) \cite{trawny2005indirect, han2011novel, ko2016sine, jing2017attitude}.  It predicts the vehicle's attitude with the gyroscope model/measurement and updates this prediction with outputs from the accelerometer and magnetometer measurement. This estimation is very accurate and widely used in practical scenarios, particularly on open-source autopilot softwares like Ardupilot and PX4. However, the Extended Kalman filter has a few limitations \cite{simon2006optimal}. Determining Kalman gain after every time interval requires two steps: propagation and update, thus requiring more computations to propagate mean and covariance, and more memory to store the results. Moreover, the EKF framework also assumes Gaussian uncertainty model, which is reasonable for uncertainty propagation over a short interval of time, requiring the EKF algorithm to run at a high rate resulting in higher processor utilization. These factors make it difficult to implement EKF in low power microprocessors.

In this paper, we propose an extended $\mathcal{H}_2$ optimal estimator that can overcome the aforementioned limitations of the extended Kalman filter. The computation involved in the proposed extended $\mathcal{H}_2$ optimal estimator requires solution of the filter state dynamics, which is of the same complexity as the mean propagation in EKF, and the estimates are obtained via simple matrix-vector multiplication. In the proposed framework, the gain is solved offline about the nominal point, which eliminates the need to solve the associated Riccati equation or the convex optimization in real-time, at the cost of reduced performance. With this compromise, we show later that the proposed extended $\mathcal{H}_2$ framework performs quite well, if not better, than the EKF with much reduced computational overhead. 

The paper is organized as follows. We first present the details of the sensor model followed by a description of conventional $\mathcal{H}_{2}$ optimal attitude estimation. This is followed by the proposed extended $\mathcal{H}_{2}$ attitude estimation algorithm.
Finally, we present simulation results using the proposed estimation and compare it with the popularly used EKF-based attitude estimation. We conclude the paper with final remarks and future research directions.

\section{Sensor Modeling}
\subsection{Gyroscope Model}
The attitude matrix in terms of an Euler angle sequence is well known \cite{schaub2005analytical} is
\begin{align}
    \begin{pmatrix} \Dot{\phi}\\ \Dot{\theta}\\ \Dot{\psi} \end{pmatrix} = \vo{T} \omegab,
\end{align}
where
\begin{align}\label{kinematic_matrix}
\vo{T}(\phi,\theta,\psi) :=
  \begin{bmatrix}
   1 & \tan\theta \sin\phi & \tan\theta \cos\phi \\
   0 & \cos\phi            & -\sin\phi \\
   0 & \sin\phi \sec\theta & \cos\phi \sec\theta
   \end{bmatrix},
\end{align}
and $\omegab$ is the angular velocity of the body with respect to inertial frame.

Gyroscope sensor measurement is developed in \cite{allen1993performance} as
\begin{subequations}
\begin{align}
& ^B \omegab = \omegab_m - \vo{b} - \vo{n}_\omegab,\\
& \dot{\vo{b}} = \vo{n_b},
\end{align}\label{gyro_m}
\end{subequations}
where $^B \omegab$ is the true angular velocity of the body frame, $\omegab_m$ is angular velocity measurement, $\vo{b}$ is the bias of gyroscope, $\vo{n}_\omegab$ is gyroscope sensor noise, and $\vo{n_b}$ represents gyroscope bias noise. In this paper we assume the gyroscope bias is non-static and model it as a random walk process.

The gyroscope model with Euler angles is then given by
\begin{multline}\label{Gyro_M}
\begin{pmatrix} \dot{\vo{\Phi}} \\ \dot{\vo{b}} \end{pmatrix}
 = \begin{bmatrix} \vo{0} & -\vo{T}(\vo{\Phi})\\ \vo{0}_{3\times3} & \vo{0}_{3\times3}\end{bmatrix}\begin{pmatrix}\vo{\Phi}\\\vo{b}\end{pmatrix} + \\
 \begin{bmatrix}-\vo{T}(\vo{\Phi}) & \vo{0}_{3\times3} \\ \vo{0}_{3\times3} & I\end{bmatrix}\begin{pmatrix}\vo{n}_\omegab\\\vo{n_b}\end{pmatrix} + \\
 \begin{pmatrix}\vo{T}(\vo{\Phi})\\\vo{0}_{3\times3}\end{pmatrix}\omegab_m,
\end{multline}
where \begin{align} \vo{\Phi}:= \begin{pmatrix} \phi \\ \theta \\ \psi \end{pmatrix}. \nonumber\end{align}

With states
\begin{align}
    \vo{x}(t) =
    \begin{bmatrix}
        \vo{\Phi}\\
        \vo{b}
    \end{bmatrix},\nonumber
\end{align}
equation (\ref{Gyro_M})  can be written as a general nonlinear dynamical system
\begin{align}
\vo{\dot{x}} = \vo{f}(\vo{x},\vo{u},\vo{w},t). \label{eqn:nldyn}
\end{align}

The noise covariance matrix for $\begin{pmatrix}\vo{n_w}\\\vo{n_b}\end{pmatrix}$ is given by
\begin{gather}
\vo{Q} =
  \begin{bmatrix}
   \vo{N}_w & \vo{0}_{3\times3}\\
   \vo{0}_{3\times3} & \vo{N}_b
   \end{bmatrix}
=
  \begin{bmatrix}
  n^2_w\vo{I}_{3\times3} & \vo{0}_{3\times3}\\
   \vo{0}_{3\times3} & n^2_b \vo{I}_{3\times3}
   \end{bmatrix}.\nonumber
\end{gather}


\subsection{Accelerometer Model}
Accelerometer sensor measurement model\cite{michel2015comparative} can be formulated as:
\begin{align}
{}^B{\vo{a}} &= \vo{a}_m - \vo{n}_a
\end{align}
with true acceleration of the body frame $^B \vo{a}$, acceleration measurement $\vo{a}_m$, and accelerometer sensor noise $\vo{n}_a$. The relationship between the gravity vector ${}^I \vo{g}$ in the inertial frame and the acceleration vector ${}^B \vo{a}$ in body frame can be formulated as
\begin{align}\label{Acc_model}
   {}^B \vo{a} = {\vo{R}^{B}_I}_{acc} (\vo{\Phi}) \hspace{0.1cm} {}^I\vo{g}
\end{align}
where $\vo{R}^B_I$ is the DCM (Direction Cosine Matrix) with 3-2-1 sequence from inertial frame($I$) to body frame($B$) as
\begin{align}
\vo{R}^B_I =&
\begin{bmatrix}\nonumber
1 &0 &0\\ 0 &\cos(\phi)& \sin(\phi)\\ 0& -\sin(\phi)& \cos(\phi)
\end{bmatrix} \times
\begin{bmatrix}
\cos(\theta)& 0& -\sin(\theta)\\ 0& 1& 0\\ \sin(\theta)& 0& \cos(\theta)
\end{bmatrix}\\ \times
&
\begin{bmatrix}
\cos(\psi)& \sin(\psi)& 0\\ -\sin(\psi)& \cos(\psi)& 0\\ 0& 0& 1
\end{bmatrix}.\nonumber
\end{align}
\subsection{Magnetometer Model}
Magnetometer sensor measurement model \cite{michel2015comparative} can be formulated  as
\begin{align}
^B{\vo{m}} &= \vo{m}_m - \vo{n}_m  \label{phi_dot}
\end{align}
with true  magnetic force  of  the  body  frame $^B \vo{m}$,  magnetometer measurement $\vo{m}_m$ and magnetometer sensor noise $\vo{n}_m$. The relationship between the Earth's magnetic vector ${}^I h$ and the local magnetic vector ${}^B \vo{m}$ can be expressed as
\begin{align}
   {}^B \vo{m} = {\vo{R}^B_I}_{mag} (\vo{\Phi}) \hspace{0.1cm} {}^I \vo{h} \label{mag_model}.
\end{align}
The noise covariance matrix for $\begin{pmatrix}\vo{n_a}\\\vo{n_m}\end{pmatrix}$ is given by
\begin{gather}
\vo{R} =
  \begin{bmatrix}
   \vo{N}_a & \vo{0}_{3\times3}\\
   \vo{0}_{3\times3} & \vo{N}_m
   \end{bmatrix}
=
  \begin{bmatrix}
  n^2_a \vo{I}_{3\times3} & \vo{0}_{3\times3}\\
   \vo{0}_{3\times3} & n^2_m \vo{I}_{3\times3}
   \end{bmatrix}. \nonumber
\end{gather}


\section{$\mathcal{H}_2$ Optimal Estimation}
We next present very briefly, the necessary background for $\mathcal{H}_2$ optimal estimation method for linear systems. We consider the following linear system,
\begin{subequations}\label{system_CT}
\begin{align}
\Dot{\vo{x}}(t) &=  \vo{A} \vo{x}(t) + \vo{B}_u \vo{u}(t) + \vo{B}_w \vo{w}(t),\\
\vo{y}(t) &= \vo{C}_y \vo{x}(t) + \vo{D}_u \vo{u}(t)+ \vo{D}_w \vo{w}(t),\\
\vo{z}(t) &= \vo{C}_z \vo{x}(t),
\end{align}
\end{subequations}
where $\vo{x} \in \mathbb{R}^{n}$, $\vo{y} \in \mathbb{R}^{l}$, $\vo{z} \in \mathbb{R}^{m}$ are respectively the state vector, the measured output vector, and the output vector of interest. Variables $\vo{w} \in \mathbb{R}^{p}$ and $\vo{u} \in \mathbb{R}^{r}$ are the disturbance and the control vectors, respectively.

With the above defined system, the $\mathcal{H}_2$ state estimator has the following form,
\begin{subequations} \label{estimation_CT}
\begin{align}
    \Dot{\hat{\vo{x}}}(t) &= \vo{A} \hat{\vo{x}}(t) + \vo{B}_u \vo{u}(t) + \vo{L} (\vo{C}_y \hat{\vo{x}}(t) + \vo{D}_u \vo{u}(t) - \vo{y}(t)),\\
    \hat{\vo{z}}(t) &= \vo{C}_z \hat{\vo{x}}(t),
\end{align}
\end{subequations}
where $\hat{\vo{x}}$ is the state estimate, $\vo{L}$ is the estimator gain, and $\hat{\vo{z}}(t)$ is the estimate of the output of interest. The error equations are then given by
\begin{subequations} \label{estimator_Linear}
    \begin{align}
    \hat{\Dot{\vo{e}}}(t) &= (\vo{A}+\vo{L} \vo{C}_y) \hat{\vo{e}}(t) + (\vo{B}_w+\vo{L} \vo{D}_w)\vo{w}(t), \textrm{and}\\
    \Tilde{\vo{z}}(t) &= \vo{C}_z \hat{\vo{e}}(t).
    \end{align}
\end{subequations}

The problem of $\mathcal{H}_2$ state estimation designs is then, given a system (\ref{estimator_Linear}) and a positive scalar $\gamma$, find a matrix $\vo{L}$ such that
\begin{align}\label{min1}
\|\vo{G}_{\Tilde{z}w}(s)\|_{2} < \gamma
\end{align}
where the transfer function $\vo{G}_{\Tilde{z}w}(s)$ of the system is:
\begin{align}\label{errorTF}
\vo{G}_{\Tilde{z} w}(s) = \vo{C}_z [s \vo{I} -( \vo{A} + \vo{L} \vo{C}_y)]^{-1} (\vo{B}_w + \vo{L} \vo{D}_w) .
\end{align}

The optimization formulation to obtain $\vo{L}$ is given by:\\[2mm]
\textbf{Theorem ($\mathcal{H}_2$ Optimal Estimation)} \cite{duan2013lmis, apkarian2001continuous} : The following two statements are equivalent:
\begin{enumerate}
    \item A solution $\vo{L}$ to the $\mathcal{H}_2$ state estimator exists.
    \item $\exists$ a matrix $\vo{W}$, a symmetric matrix $\vo{Q}$, and a symmetric matrix $\vo{X}$ such that:
\end{enumerate}
        \begin{align}\label{CT_LMI}
        \begin{bmatrix}\nonumber
        \vo{XA}+\vo{W} \vo{C}_y +(\vo{XA}+\vo{W} \vo{C}_y)^T &  \vo{X} \vo{B}_w+\vo{W} \vo{D}_w\\
        \vo{*}      & -\vo{I}
        \end{bmatrix}
        < 0,\\
        \begin{bmatrix}\nonumber
        \vo{-Q}  &  \vo{C}_z\\
        \vo{*} & \vo{-X}
        \end{bmatrix}
        < 0,\\
         \textbf{trace}(\vo{Q}) < \gamma^2.
        \end{align}

The $\mathcal{H}_2$ optimal estimator gain is recovered by $\vo{L} = \vo{X}^{-1} \vo{W}$. This optimal gain ensures that:
\begin{align}
\vo{e}(t) = \vo{x}(t) - \hat{\vo{x}} \rightarrow 0 \hspace{0.5cm} \text{ as } t \rightarrow \infty,
\end{align}
 and ensures that $\hat{\vo{x}}(t)$ is an asymptotic estimate of $\vo{x}(t)$.
\\
\section{Extended $\mathcal{H}_2$ Estimation}
\begin{figure}[h!]
\centering
\includegraphics[trim={1cm 0.5cm 1cm 0cm},scale=0.45]{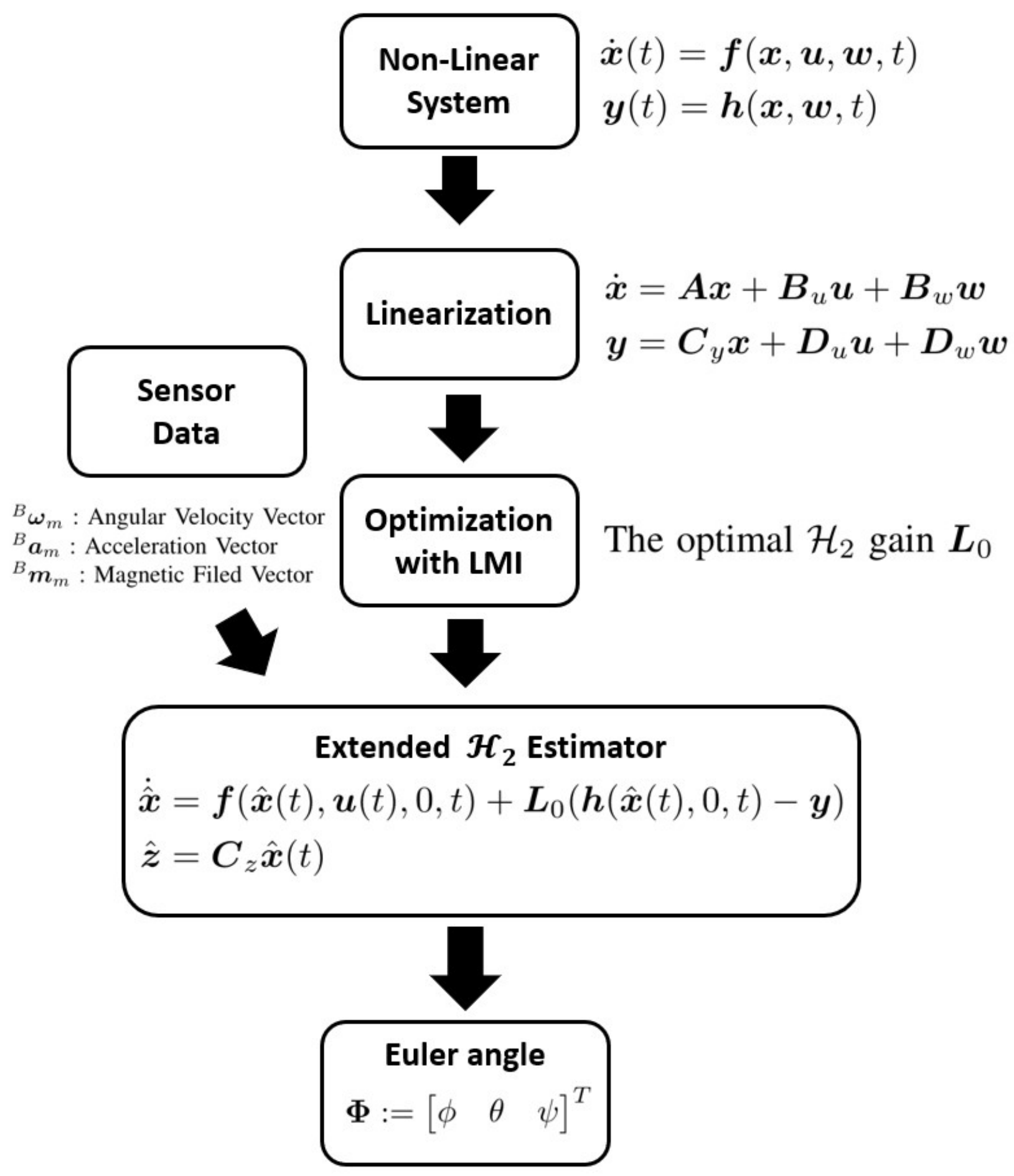}
\caption{Estimation algorithm}
\label{algorithm}
\end{figure}

The proposed extended $\mathcal{H}_2$ estimation framework is summarized in Fig. \ref{algorithm}. Recall from (\ref{eqn:nldyn}) that the gyro model is nonlinear, and can be represented as
\begin{align}\label{sys_eq}
\Dot{\vo{x}} &= \vo{f}(\vo{x},\vo{u},\vo{w},t),
\end{align}

where $\vo{u}(t) := \omegab_m(t)$ and $\vo{w}(t):=\begin{bmatrix}
        \vo{n}_\omegab(t)\\
        \vo{n}_b(t)
    \end{bmatrix}$.

The measurement equations with accelerometer and magnetometer model (\ref{Acc_model}, \ref{mag_model}) can be written as the following nonlinear output equation
\begin{align} \label{eqn:measureEQ}
\vo{y}(t) &= \vo{h}(\vo{x},\vo{w},t) = \vo{R}^B_I(\vo{\Phi})
    \begin{bmatrix}
       \vo{g} \\ \vo{h}
   \end{bmatrix}  + \vo{D}_w \vo{w}(t),
\end{align}
where,
\vspace{-.3cm}
\begin{gather*}
\vo{R}^B_I(\vo{\Phi})
   =
    \begin{bmatrix}
       \vo{R}_{acc}(\vo{\Phi})& \vo{R}_{mag}(\vo{\Phi})
    \end{bmatrix},
\vo{D_w} =
   \begin{bmatrix}
        \vo{I}_{3\times3} & \vo{0}_{3\times3}\\
        \vo{0}_{3\times3} & \vo{I}_{3\times3}
  \end{bmatrix}.
\end{gather*}

Extended $\mathcal{H}_2$ estimation is simply $\mathcal{H}_2$ estimation extended to nonlinear system models, along the lines of the extended Kalman filtering. In extended Kalman filtering, the uncertainty is propagated using the linear system along the state trajectory, and the Kalman gain is computed at every time step with the instantaneous linear system. In the extended $\mathcal{H}_2$ framework, in theory, we can solve for the optimal $\mathcal{H}_2$ gain along the trajectory, however this may be computationally prohibitive for cheap processors. Instead, we design
the optimal $\mathcal{H}_2$ gain about the \textit{nominal} operating point, but use the \textit{nonlinear system dynamics} to evolve the estimator's states.


A linear approximation is implemented at nominal operating point $(\vo{x}_0,\vo{u}_0, \vo{w}_0) = \vo{0}$. Linearizing about this nominal point, we get
\begin{align}
    \vo{f}&(\vo{x}(t), \vo{u}(t), \vo{w}(t), t) \approx \vo{f}(\vo{x}_0, \vo{u}_0,\vo{w}_0, t)\\ \nonumber
    &+\frac{\partial \vo{f} (\vo{x},\vo{u},t)}{\partial \vo{x}}\Bigg|_{\text{nominal}} \vo{x}(t)  +\frac{\partial \vo{f} (\vo{x},\vo{u},t)}{\partial \vo{u}}\Bigg|_{\text{nominal}} \vo{u}(t)\\ \nonumber
    &+ \frac{\partial \vo{f} (\vo{x},\vo{u},t)}{\partial \vo{w}}\Bigg|_{\text{nominal}} \vo{w}(t) + \textrm{H.O.T.}
\end{align}
Jacobian matrix of the system is defined as:
\begin{align} \nonumber
\vo{A}  := \frac{\partial \vo{f}(\vo{x},\vo{u},t)}{\partial \vo{x}}\Bigg|_{\text{nominal}},
\vo{B}_u  := \frac{\partial \vo{f}(\vo{x},\vo{u},t)}{\partial \vo{u}}\Bigg|_{\text{nominal}},\\
\vo{B}_w  := \frac{\partial \vo{f}(\vo{x},\vo{u},t)}{\partial \vo{w}}\Bigg|_{\text{nominal}}.
\end{align}
Then the linearization of equation (\ref{sys_eq}) is
 \begin{align}
     \Dot{\vo{x}}(t) &= \vo{f}(\vo{x}(t),\vo{u}(t),\vo{w}(t),t)\\
           &\approx \vo{f}(\vo{x_0},\vo{u_0}, t) + \vo{A} \vo{x}(t) + \vo{B_u} \vo{u}(t) + \vo{B_w} \vo{w}(t) .\nonumber
 \end{align}

The measurement model \eqref{eqn:measureEQ} can be similarly approximated as:
\begin{align}
    \vo{h}(\vo{x}(t), \vo{w}(t), t) \approx \nonumber\\
    \vo{h}(\vo{x_0}, \vo{w_0}, t) &+\frac{\partial \vo{h} (\vo{x},t)}{\partial \vo{x}}\Bigg|_{\text{nominal}} \vo{x}(t)\nonumber\\
    &+ \frac{\partial \vo{h} (\vo{x},t)}{\partial \vo{v}}\Bigg|_{\text{nominal}} \vo{w}(t) + \textrm{H.O.T.}
\end{align}
Jacobian matrix of the measurement is defined as:
\begin{align}
\vo{C}_y  := \frac{\partial \vo{h} (\vo{x},t)}{\partial \vo{x}}\Bigg|_{\text{nominal}},
\vo{D}_w  := \frac{\partial \vo{h} (\vo{x},t)}{\partial \vo{w}}\Bigg|_{\text{nominal}}.\nonumber
\end{align}
The linearized measurement equation \eqref{eqn:measureEQ} can then be written as:
 \begin{align}
     \vo{y}(t) &= \vo{h}(\vo{x}(t),\vo{w}(t),t)\\
           &\approx \vo{h}(\vo{x_0},\vo{w_0}, t) + \vo{C}_y \vo{x}(t) + \vo{D}_w(t) \vo{w}(t).\nonumber
 \end{align}

The linear system, about the nominal operating point, is therefore
\begin{subequations}\label{Lin_systems}
\begin{align}
\Dot{\vo{x}} &=  \vo{A}  \vo{x} + \vo{B}_u  \vo{u} + \vo{B}_w  \vo{w},\\
\vo{y} &= \vo{C}_y  \vo{x} + \vo{D}_u  \vo{u}+ \vo{D}_w  \vo{w}.
\end{align}
\end{subequations}
The optimal $\mathcal{H}_2$ gain $\vo{L}_0$ can then be determined by solving the optimization problem in (\ref{CT_LMI}), where the subscript $0$ is used to indicate that it is determined about the nominal operating point.

Once the gain $\vo{L}_0$ is determined it is used to implement the $\mathcal{H}_2$ filter for the nonlinear system. We present a new implementation, called the extended $\mathcal{H}_2$ filter, where the filter states are propagated using the nonlinear dynamics. In conventional $\mathcal{H}_2$ filters, the error propagation occurs using the linear system. The filter dynamics and output equation for the extended $\mathcal{H}_2$ filter are given by
\begin{subequations}\label{estimator}
    \begin{align}
    \Dot{\Hat{\vo{x}}} &= \vo{f}(\hat{\vo{x}}(t),\vo{u}(t),0,t) +  \vo{L}_0(\vo{h}(\Hat{\vo{x}}(t),0,t) - \vo{y}),\\
    \hat{\vo{z}} &= \vo{C}_z \hat{\vo{x}}(t).
    \end{align}
\end{subequations}



\section{Result}

\subsection{Simulation Set Up}
\begin{figure}[h!]
\centering
\includegraphics[trim={0 .5cm 0 .5cm},width=0.47\textwidth]{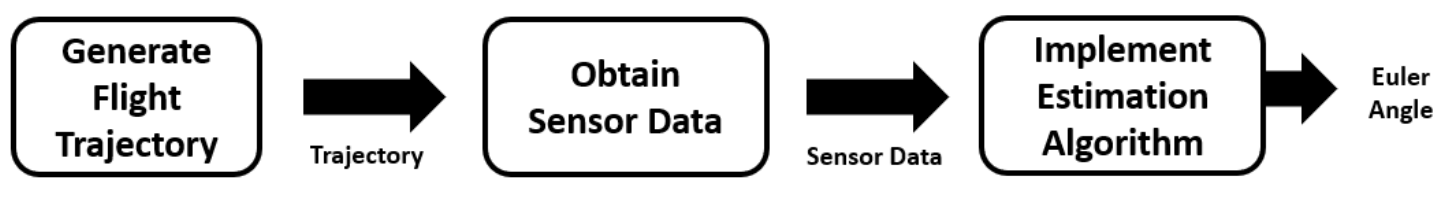}
\caption{Simulation Flow Chart}
\label{sim_chart}
\end{figure}

The proposed  extended $\mathcal{H}_2$ filter is applied to the attitude estimation problem, and its performance compared with extended Kalman filter based attitude estimation. The comparison is done in terms of estimation accuracy and computational time in a MATLAB based simulation environment, as shown in Fig. \ref{sim_chart}.

\begin{figure}[h!]
\centering
\includegraphics[width=0.47\textwidth]{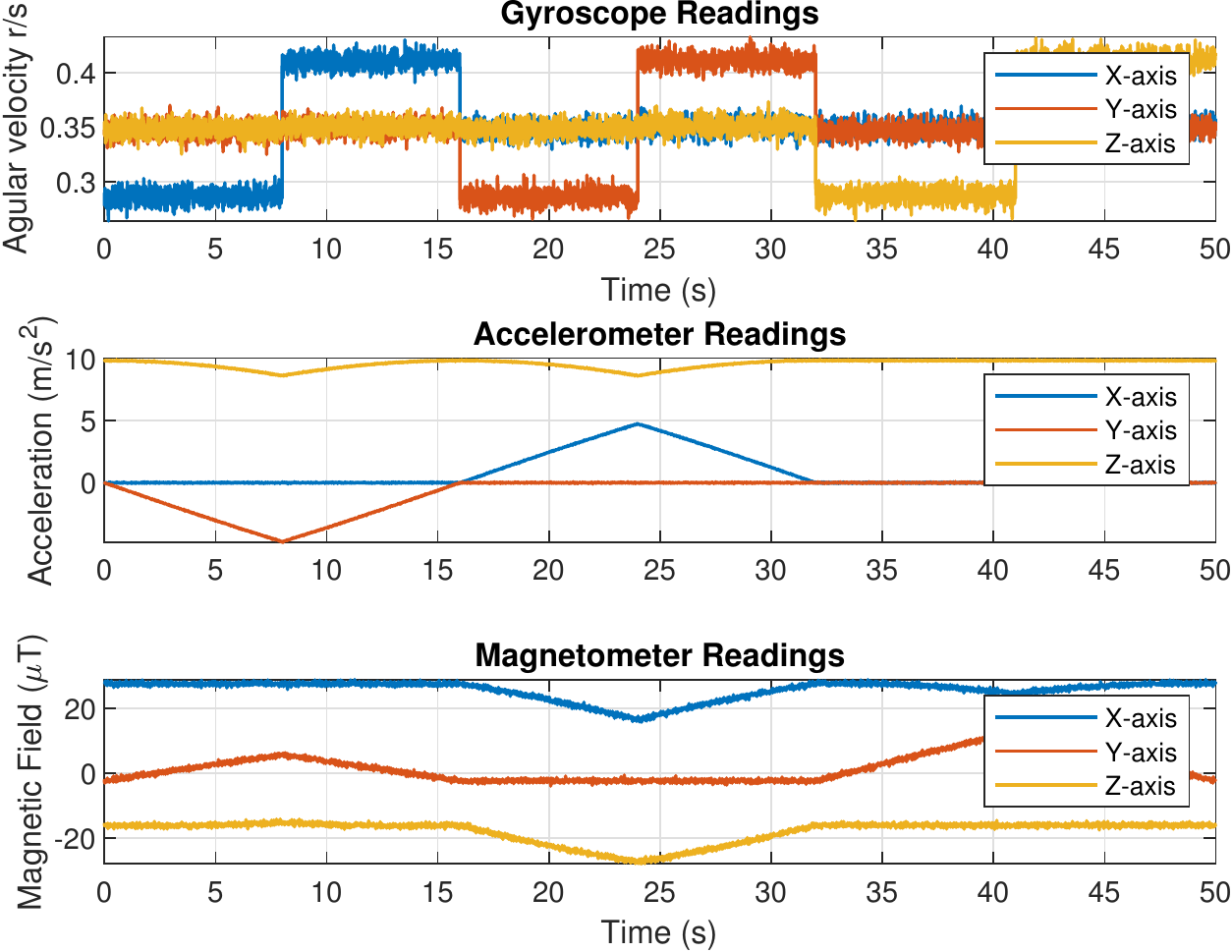}
\caption{IMU Data from MATLAB function \textsl{imuSensor} with MPU- 9250 data}
\label{IMU_Low}
\end{figure}

Sensor characteristics like noise level, bias, etc. are set by sensor data sheet \cite{MPU} of MPU 9250, an affordable commercial sensor. Moreover, this sensor is used in the Pixhawk flight computer, which is broadly used in unmanned aerial systems. In the MATLAB simulation, IMU data is generated by the \texttt{imuSensor} function \cite{MAT}. The data is shown in Fig. \ref{IMU_Low}. Two flight scenarios are used to verify the proposed extended $\mathcal{H}_2$ estimation algorithm.\\

\begin{figure}[h!]
\centering
\includegraphics[width=0.47\textwidth]{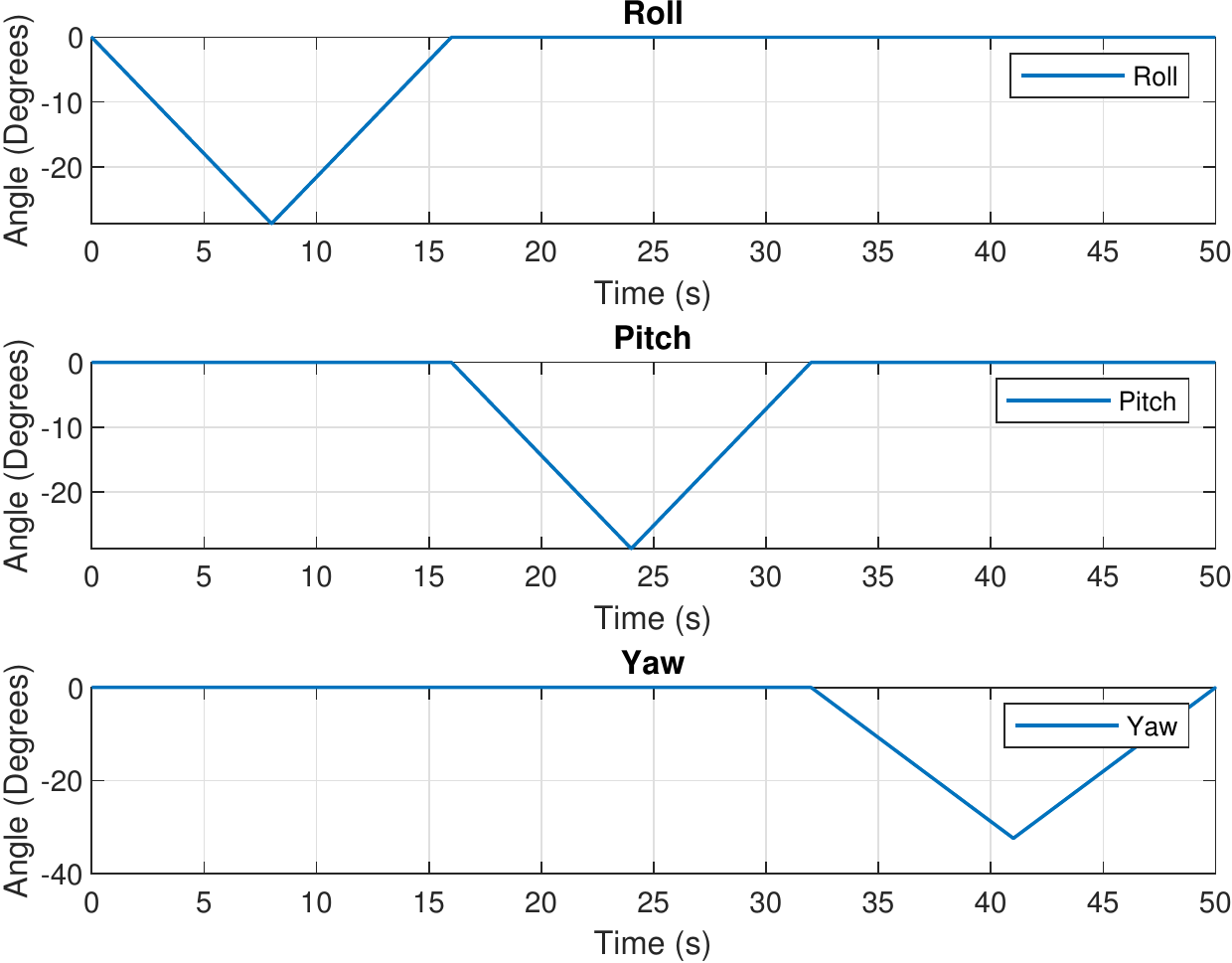}
\caption{True trajectories of Euler angles for Case I.}
\label{Tra_Low}
\end{figure}

\textbf{Case I: Slow and Small Angle Movements} -- Here we consider angular movements $< 30^{\circ}$ about all three axes of the vehicle independently. This case broadly covers forward/backward and left/right cruise flight of popular quad-rotor based UAVs. Simulation is run for a time duration of 50 seconds and with an angular velocity $\pi/50$ rad/s. The simulated true state trajectories are shown in Fig. \ref{Tra_Low}.

\begin{figure}[h!]
\centering
\includegraphics[width=0.47\textwidth]{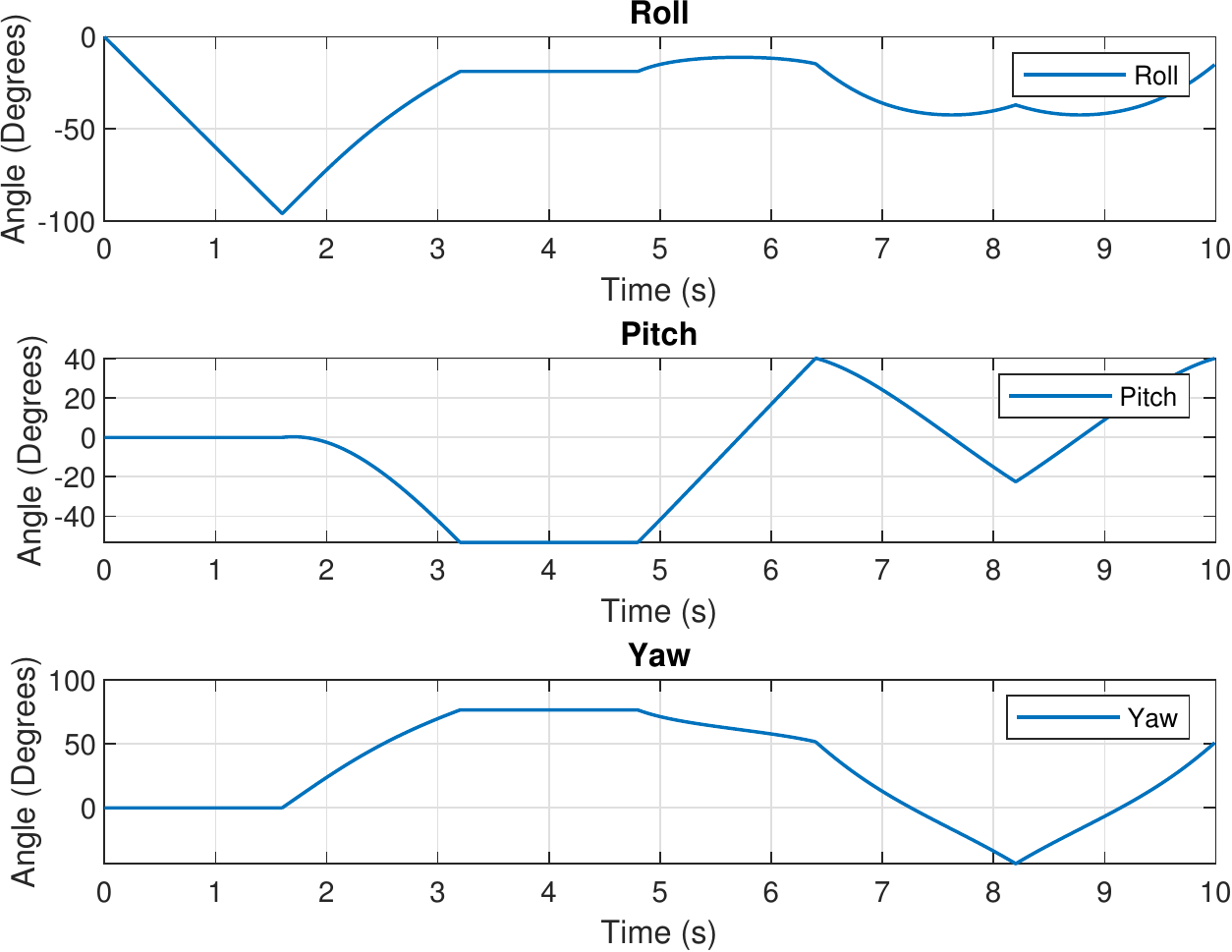}
\caption{True trajectories of Euler angles for Case II.}
\label{Tra_high}
\end{figure}

\newpage
\textbf{Case II: Fast and Large Angle Movements} -- Here we consider angular variation $> 30^{\circ}$  about all three axes of the vehicle simultaneously. It represents scenarios of rapid movements or motion in the presence of wind disturbance during flight or aggressive maneuvers. Simulation is run for a duration of 10 seconds and with an angular velocity of $\pi/3$ rad/s. The simulated true state trajectories are shown in Fig. \ref{Tra_high}.

\subsection{Simulation Results}
\begin{figure}[h]
\centering
\includegraphics[trim={1cm 1cm 1cm .6cm},width=0.47\textwidth]{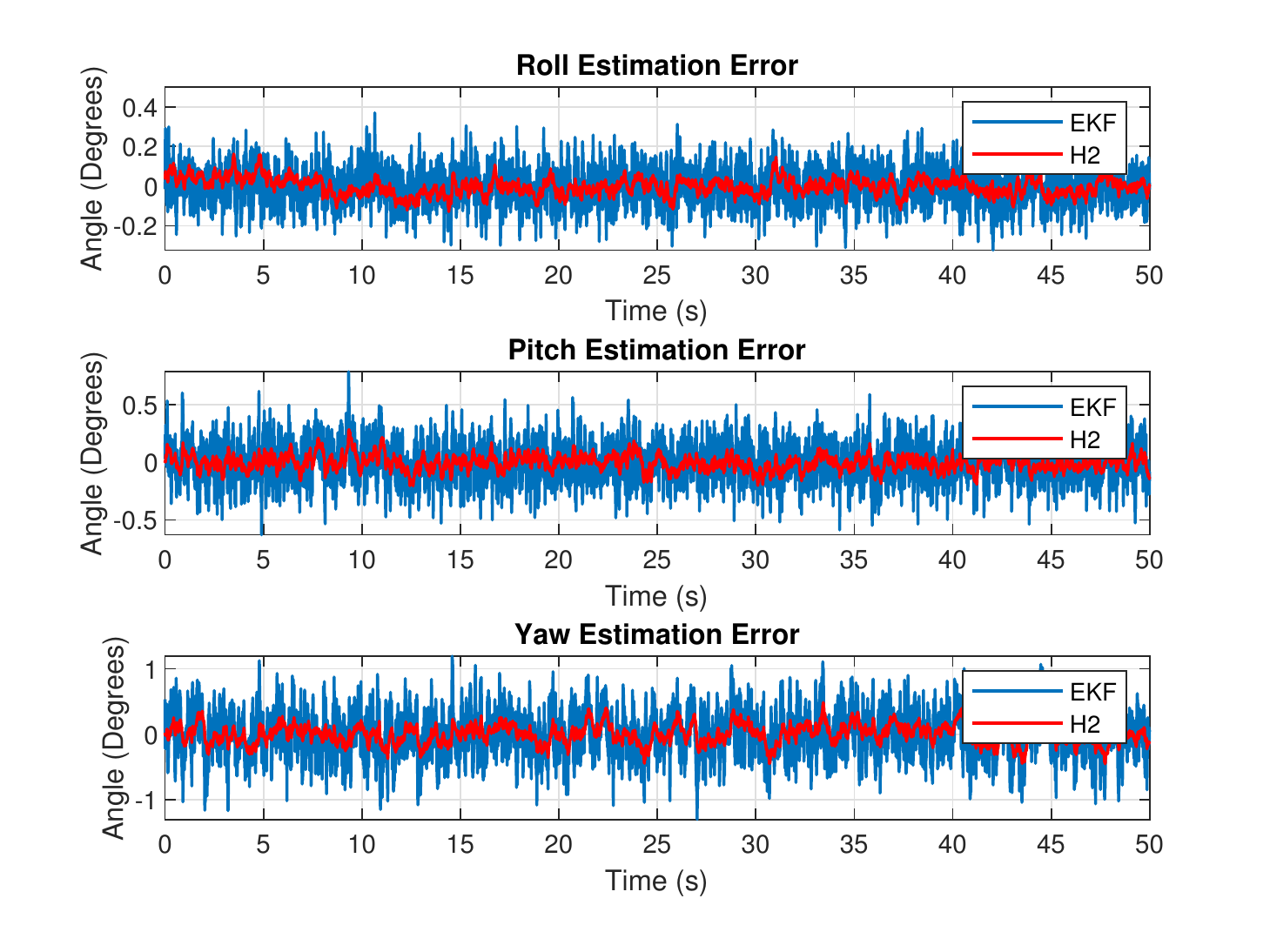}
\caption{Comparison of the proposed extended $\mathcal{H}_2$ filter and the EKF for Case I.}
\label{comp_low}
\end{figure}

Let us examine the performance of the extended $\mathcal{H}_2$ estimator with that of the standard EKF in terms of root mean squared (RMS) error, memory use, and computational time required,

\textbf{Case I:} The comparison of the two estimators for Case I is shown in Fig. \ref{comp_low}. We observe that the error of extended $\mathcal{H}_2$ estimator is less than that of EKF. The RMS error, and the upper and lower bounds of the error, for both the filters are shown in table \ref{RMS_low}. From the plots and the data in the tables, we can conclude that the performance of extended $\mathcal{H}_2$ estimator is better than that of EKF.

\begin{table}[h]
\caption{RMS error for Case I.} \label{RMS_low}
\vspace{-0.3cm}
\begin{center}
\renewcommand{\arraystretch}{1.5}
\begin{tabular}{|c||c|c|c|}
\hline
Algorithm & Roll angle (${}^\circ$) & Pitch angle (${}^\circ$) & Yaw angle (${}^\circ$)\\
\hline \hline
Extended $\mathcal{H}_2$ & 0.0331 & 0.0538 & 0.1107\\
\hline
EKF & 0.0533 & 0.0988 & 0.2298\\
\hline
\end{tabular}
\end{center}
\end{table}

\begin{table}[h]
\caption{Min-max error for Case I.} \label{min-max_low}
\vspace{-0.2cm}
\begin{center}
\renewcommand{\arraystretch}{1.5}
\begin{tabular}{|c||c|c|c|}
\hline
Algorithm & Roll angle (${}^\circ$) & Pitch angle (${}^\circ$) & Yaw angle (${}^\circ$)\\
\hline \hline
Extended $\mathcal{H}_2$ &  [-0.152 0.125] & [-0.284 0.240] & [-0.423 0.402]\\
\hline
EKF & [-0.257 0.249] & [-0.523 0.444]& [-1.087 0.937]\\
\hline
\end{tabular}
\end{center}
\end{table}

\begin{figure}[h]
\centering
\includegraphics[trim={1cm 1cm 1cm .5cm},width=0.47\textwidth]{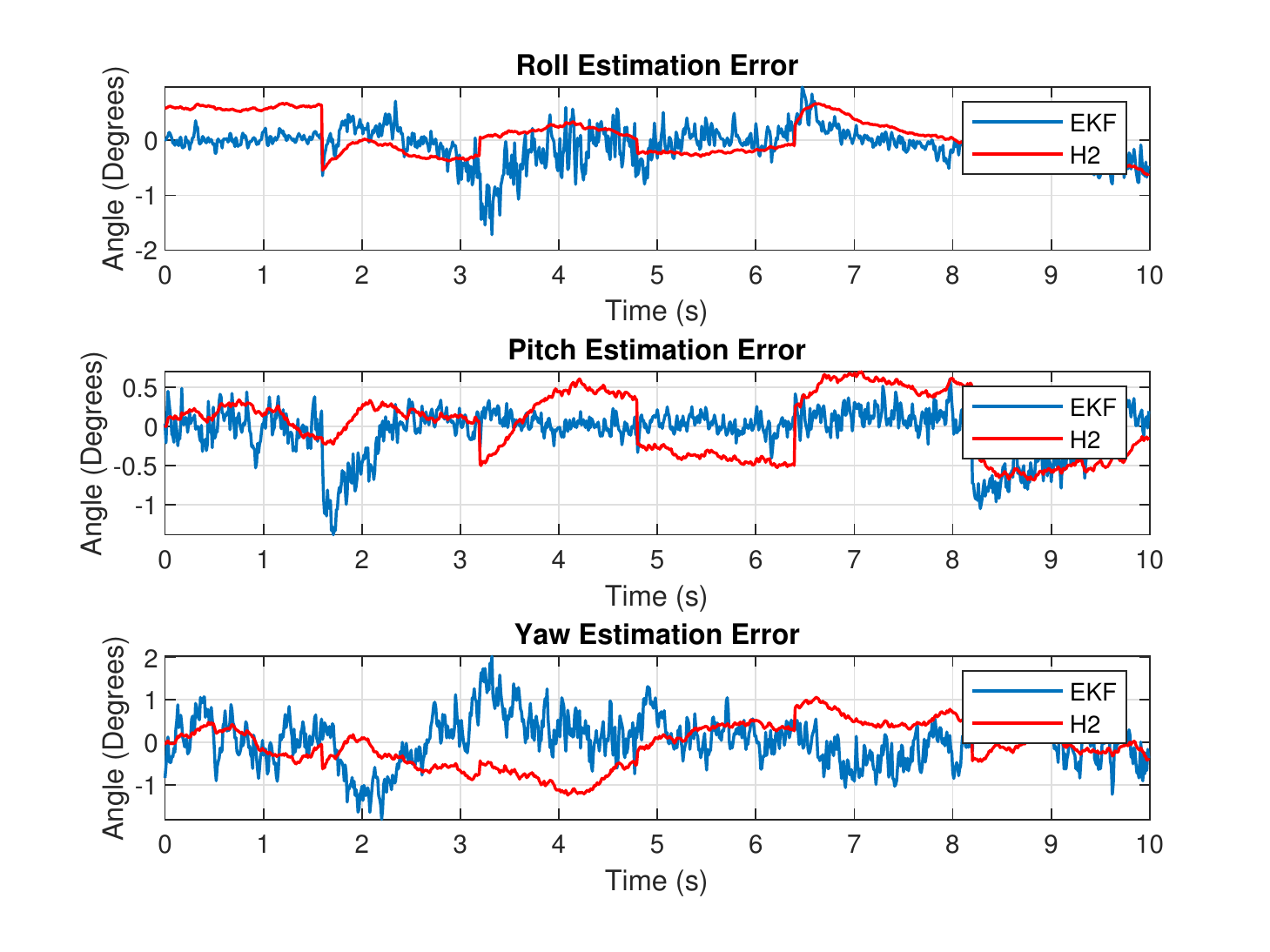}
\caption{Comparison of the proposed extended $\mathcal{H}_2$ filter and the EKF for Case II.}
\label{comp_high}
\end{figure}

\textbf{Case II:} The results of the two estimators for case II are shown in Fig. \ref{comp_high}. From table \ref{min-max_high}, we observe that the extended $\mathcal{H}_2$ filter has lower error bounds. From table \ref{RMS_high}, we observe that the RMS errors are comparable for both the filters.

\vspace{-0.4cm}

\begin{table}[h]
\caption{RMS error for Case II.} \label{RMS_high}
\vspace{-0.3cm}
\begin{center}
\renewcommand{\arraystretch}{1.5}
\begin{tabular}{|c||c|c|c|}
\hline
Algorithm & Roll angle (${}^\circ$) & Pitch angle (${}^\circ$) & Yaw angle (${}^\circ$)\\
\hline \hline
Extended $\mathcal{H}_2$ & 0.3045 & 0.3260 & 0.3121\\
\hline
EKF & 0.3073 & 0.2656 & 0.5275\\
\hline
\end{tabular}
\end{center}
\end{table}

\begin{table}[h!]
\caption{Min-max error for Case II.} \label{min-max_high}
\vspace{-0.3cm}
\begin{center}
\renewcommand{\arraystretch}{1.5}
\begin{tabular}{|c||c|c|c|}
\hline
Algorithm & Roll angle (${}^\circ$) & Pitch angle (${}^\circ$) & Yaw angle (${}^\circ$)\\
\hline \hline
Extended $\mathcal{H}_2$ &  [-0.592    0.730] & [-0.897    0.602] & [-0.857 0.664]\\
\hline
EKF & [-2.212   1.404] & [-1.84000.692] & [-2.227    3.097]\\
\hline
\end{tabular}
\end{center}
\end{table}

\newpage
Based on the estimation errors for Case I and Case II, we can conclude that both the filters perform equally well, with the extended $\mathcal{H}_2$ filter performing slightly better. The main advantage of the  extended $\mathcal{H}_2$ filter is in the implementation efficiency. The results for the average execution time are shown in table \ref{Com_time}, which reveals that the extended $\mathcal{H}_2$ estimator requires $50\%$ less computational time than EKF, making it $2\times$ more efficient than EKF. Table \ref{Com_time} also shows the variability in the computational time, which is about $3\times$ less. The reduced variability in the computational time increases the reliability of the real-time tasks that will execute in the microprocessor. Typically, more time is allotted to real-time tasks with large variability in computational time. This further improves the computational efficiency of the proposed extended $\mathcal{H}_2$ estimator.

\begin{table}[h!]
\caption{Computational Time Comparison} \label{Com_time}
\vspace{-0.3cm}
\begin{center}
\renewcommand{\arraystretch}{1.5}
\setlength{\tabcolsep}{10pt}
\begin{tabular}{|c||c|c|}
\hline
Algorithm & Mean Time ($ms$) & Standard Deviation ($ms$)\\
\hline \hline
Extended $\mathcal{H}_2$ &  0.853 & 0.244\\
\hline
EKF & 1.7 & 0.736 \\
\hline
\end{tabular}
\end{center}
\end{table}

\section{CONCLUSIONS}
This paper presents a new nonlinear estimation framework, based on $\mathcal{H}_2$ optimal state estimation, for attitude estimation in low power microprocessors. We showed that the performance of the proposed estimator is comparable, if not better, than that of the EKF algorithm which is popularly used in the application space considered. The primary advantage of the proposed framework is the $2\times$ computational efficiency, and the $3\times$ robustness with respect to computational uncertainty. Both these factors make the proposed attitude estimation algorithm very attractive for small UAVs with low power microprocessors. The present work is based on Euler angles for attitude estimation, which is known to have singularities. In our future work, we will address attitude estimation using quaternions in the extended $\mathcal{H}_2$ framework.

\addtolength{\textheight}{-12cm}   


\bibliographystyle{elsarticle-num}
\bibliography{H2filtering}

\end{document}